\documentclass[12pt]{article}%
\usepackage[applemac]{inputenc}
\usepackage{amsmath,amssymb}
\usepackage{amsfonts}
\usepackage{amsmath,amssymb}
\usepackage[applemac]{inputenc}
\usepackage{color}
\usepackage{color, colortbl}
\usepackage{amsmath}
\usepackage{amssymb}
\usepackage{bbm}
\usepackage{graphicx}  
\usepackage{subcaption}  
\usepackage{tikz}
\usepackage{geometry}
\usetikzlibrary{arrows}
\usepackage{amsfonts}
\usepackage{amsmath,amssymb}
\usepackage[applemac]{inputenc}
\usepackage{color}
\usepackage{amsmath}
\usepackage{amssymb}
\usepackage{graphicx}
\usepackage{tikz}
\newtheorem{theorem}{Theorem}[section]
\newtheorem{lemma}[theorem]{Lemma}
\newcommand{\E}{E}

\newtheorem{definition}[theorem]{Definition}

\newtheorem{remark}[theorem]{Remark}

\usetikzlibrary{arrows,shapes,automata,petri}

\newcommand{\dproof}{\noindent {Proof.} \quad}
\newcommand{\fproof}{\hfill $\square$ \bigskip}

\numberwithin{equation}{section}

\definecolor{LightCyan}{rgb}{0.88,1,1}

\def\1B{\text{1\!\!I}}

\def\<{\langle}
\def\>{\rangle}

\def\E{\mathbb{E}}
\def\A{\mathcal{A}}
\def\R{\mathbb{R}}

\begin{document}
\title{SPDE Games Driven by a Brownian Sheet with Applications to Pollution Minimization}
\author{Nacira Agram$^{1},$ Bernt \O ksendal$^{2}$, Frank Proske$^{2}$, and Olena Tymoshenko$^{2,3}$}
\date{3 March 2025}
\maketitle

\footnotetext[1]{Department of Mathematics, KTH Royal Institute of Technology 100 44, Stockholm, Sweden. \newline
Email: nacira@kth.se. Work supported by the Swedish Research Council grant (2020-04697).}

\footnotetext[2]{%
Department of Mathematics, University of Oslo, Norway. \\
Emails: oksendal@math.uio.no, proske@math.uio.no, olenaty@math.uio.no}

\footnotetext[3]{%
Department of Math. An. \& Prob. Th., Igor Sikorsky Kyiv Polytechnic Institute, Ukraine.}
\begin{abstract}
This paper studies a nonzero-sum stochastic differential game in the context of shared spatial-domain pollution control. The pollution dynamics are governed by a stochastic partial differential equation (SPDE) driven by a Brownian sheet, capturing the stochastic nature of environmental fluctuations. Two players, representing different regions, aim to minimize their respective cost functionals, which balance pollution penalties with the cost of implementing control strategies.  

The nonzero-sum framework reflects the interdependent yet conflicting objectives of the players, where both cooperation and competition influence the outcomes. We derive necessary and sufficient conditions for Nash equilibrium strategies, using a maximum principle approach. This approach involves the introduction of a new pair of adjoint variables, \( (L_1, L_2) \), which do not appear in a corresponding formulation with the classical (1-parameter) Brownian motion.  

Finally, we apply our results to two case studies in pollution control, demonstrating how spatial and stochastic dynamics shape the equilibrium strategies.
\end{abstract}




\section{Introduction}

The global challenge of pollution management often involves multiple stakeholders or regions sharing a common geographical domain. Each entity seeks to mitigate pollution while balancing the economic and operational costs of their interventions. In such scenarios, the interactions between stakeholders are interdependent, as one entity's actions influence the outcomes for others. These dynamics can be modeled effectively using nonzero-sum game theory, which accommodates both competitive and cooperative behaviors among players.

In this work, we develop a stochastic differential game framework for pollution control, focusing on a two-player, nonzero-sum setting. The state dynamics are governed by a stochastic partial differential equation (SPDE) driven by  a time-space Brownian motion (Brownian sheet, see (a)-(b) visualizations), which captures spatial and temporal randomness. This approach reflects real-world uncertainties such as fluctuating weather conditions, irregular external emissions, and spatial variations in pollution levels. Each player, representing region or entity, aims to minimize the individual cost functional, which includes both  penalties  for excessive pollution and the economic  costs of control measures.

Unlike zero-sum games, where one player's gain is exactly the other's loss, nonzero-sum games allow for a broader range of strategic interactions. In the pollution control context, this framework captures the reality that both players aim to achieve mutually beneficial outcomes while accounting for their individual objectives. The Nash equilibrium serves as the solution concept, representing a state where no player can unilaterally improve their outcome by altering their strategy.

Building on prior works, such as the investigation of forward-backward stochastic differential games by \O ksendal \& Sulem (2019) \cite{OS} and the analysis of control strategies in stochastic differential games by
 \O ksendal \& Reikvam (1998) \cite{OR}, we further extend the analysis of stochastic differential games in the context of SPDEs. 
 
 Stochastic differential games have been widely applied in various domains. For instance, in recent research Zuliang et al. (2024) \cite{Zh2} investigate  the dynamics of pollution control through stochastic evolution games, emphasizing the role of cooperative strategies. Li \& Zhang (2023) \cite{LZ} present an optimal control framework for industrial pollution, utilizing stochastic differential equations to model real-world uncertainties such as equipment reliability and raw material variability.

To illustrate  a real-world applications, stochastic differential games are are frequently used in financial markets to model the strategic interactions of competing investors optimizing their portfolios under uncertain market conditions ( e.g. see Fleming \& Hernandez-Hernandez (2003) \cite{F}). In distributed control systems, the system can describe competing agents managing energy distribution in smart grids, balancing efficiency and reliability under stochastic demands  Mertikopoulos \& Staudigl (2018)\cite{MS}. 
 In the field of financial mathematics, Cardaliaguet \& Rainer (2009) \cite{CR} investigate stochastic differential games in portfolio optimization, where competing investors adjust their strategies in response to market fluctuations. In energy networks, stochastic differential games can model strategic interactions in electricity pricing, where noise represents fluctuations in energy supply and demand.  Bensoussan et al. (2015) \cite{BGQ} explore applications in energy markets, modeling competition between producers in electricity pricing using mean-field game approaches.

 Furthermore, recent results obtained by Agram et al. (2025) \cite{AOPT1}, where the authors examine optimal control of SPDEs driven by time-space Brownian motion, provide new opportunities to study stochastic differential games of SPDEs under such dynamics.
\begin{figure}[h]
    \centering
    \begin{subfigure}[t]{0.49\textwidth}
        \includegraphics[width=\linewidth]{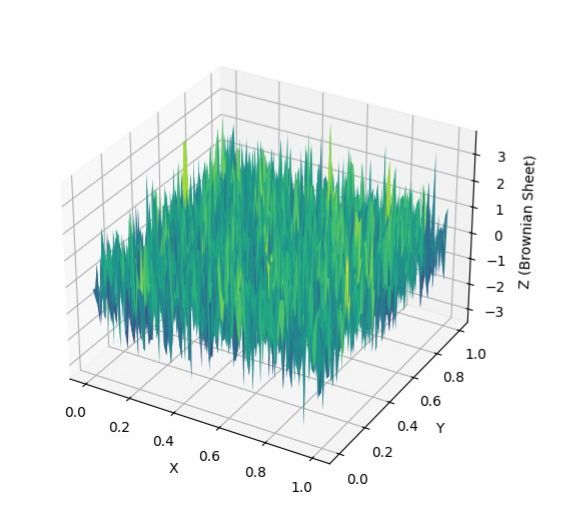}
        \caption{Visualization of Brownian Sheet}
        \label{fig:B}
    \end{subfigure}
    \hfill
    \begin{subfigure}[t]{0.44\textwidth}
        \includegraphics[width=\linewidth]{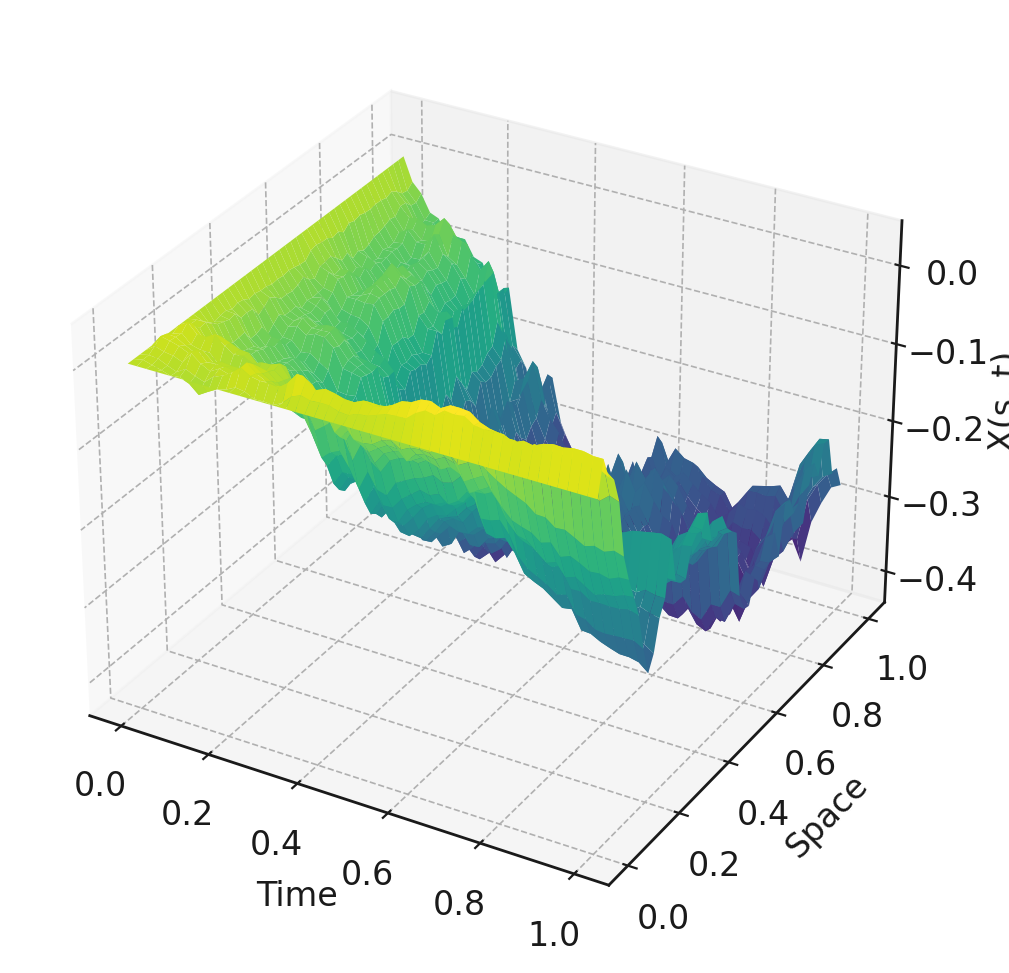}
        \caption{SPDE driven by Brownian Sheet}
    \end{subfigure}
\end{figure}

Expanding on the results of Agram et al. (2025) \cite{AOPT1}, we investigate stochastic differential games of SPDEs driven by time-space Brownian motion. In this study, we establish necessary and sufficient conditions for Nash equilibria, with the Hamiltonians of each player serving as key components in defining optimal strategies. To demonstrate the practical implications of our theoretical results, we apply them to a pollution control problem, deriving equilibrium strategies within a stochastic framework.

The paper is organized as follows:

In Section 2, we explain the basic concepts needed for this study. We introduce an SPDEs driven by time-space Brownian motion. We also describe key ideas such as Brownian sheets, filtration structures, and It\^o calculus in the plane. These concepts are important for our analysis.

In Section 3, we present the main theoretical framework for nonzero-sum stochastic differential games. We describe the pollution control problem and how the state of the system follows a SPDE. We also define the strategies of the two players and find the necessary and sufficient conditions for Nash equilibrium.
In Section 4, we apply our theoretical framework to two pollution control problems involving two competing regions. Each region seeks to minimize pollution while balancing the costs of control efforts, modeled using SPDE. We derive optimal control strategies and Nash equilibrium conditions by means of necessary and sufficient conditions. A special case with symmetric parameters is analyzed, showing how external factors and randomness influence decision-making.

\section{Preliminaries}
Let \( \zeta = (s, a) \in \mathbb{R}^2 \) denote a 2-dimensional index, where \( s \) denotes time and \( a \) represents space. The differential element corresponding to \( \zeta \) is written as \( d\zeta = ds \, da \), which signifies the 2-dimensional Lebesgue measure over the time-space domain.

The process \( B(t, x) \) represents a Brownian sheet, a generalization of Brownian motion to two dimensions. Here $B(t,x)$ is indexed by time \( t \geq 0 \) and space \( x \in \mathbb{R} \). This Brownian sheet is assumed to satisfy the following properties:
\begin{enumerate}
    \item [(i)] \textbf{continuity and independence:} for \( t \geq 0 \), \( B(t, \cdot) \) is a family of continuous functions indexed by space \( x \in \mathbb{R} \). The increments of this process over disjoint regions of $t$  are independent.
    \item [(ii)] \textbf{stationarity in space:} for each fixed \( t \), the process \( B(t, x) \) is stationary and its increments are independent  over disjoint intervals of $x$.
\end{enumerate}

To define the filtration associated with the Brownian sheet, let \( \mathcal{F}_{t,x} \)  represent  the natural filtration generated by the Brownian sheet \( B(s, a) \) up to time \( t \) and space \( x \), i.e.,
\[
\mathcal{F}_{t,x} = \sigma(B(s, a), 0 \leq s \leq t, a \in [0, x]),
\]
which contains all the information available to the point \( (t,x) \) in the time-space domain.

We consider a process \( Y(t, x) \) driven by a Brownian sheet, which is defined as:
\begin{align} \label{eq}
Y(t, x) &= Y_0 + \int_{R_{t,x}} \alpha(s,a) \, d(s,a) + \int_{R_{t,x}} \beta(s,a) \, B(d(s,a))\nonumber\\&+ \iint_{R_{t,x} \times R_{t,x}} \psi(s,a, s',a') \, B(d(s,a)) B(d(s',a')) ,
\end{align}
where \( Y_0 \) is the initial value of the process, \( \alpha(s,a) \) is the drift term, \( \beta(s,a) \) is the diffusion term, and \( \psi(s,a, s',a') \) describes the interaction term between the components of the Brownian sheet. 

The  first integral in (\ref{eq}) represents the deterministic change in $Y(t,x)$ due to $\alpha(t,x)$, the second accounts for stochastic fluctuations introduced by the Brownian  sheet,  and the third captures the correlation effects between parts of the Brownian sheet.  The integration is performed over a  region \( R_{t,x} = [0,t] \times [0,x] \). These integrals were rigorously developed in Cairoli \& Walsh (1975)\cite{CW} and Wong \& Zakai (1974) \cite{WZ74}, and they form the foundation of stochastic calculus in the plane. 

\subsection{The It\^o Formula}

To study the necessary and sufficient conditions for Nash equilibrium problems, we utilize a version of the It\^o formula adapted to such systems. To ensure clarity in notation, we adopt the framework introduced by Wang \& Zakai (1978)\cite{WZ}. For simplicity, we define \( z = (t, x) \) and \( R_z = [0, t] \times [0, x] \), which represents a rectangular region in the time-space domain.  The integral \( \int_{R_z} \phi(\zeta) B(d\zeta) \) is interpreted as the It\^o integral with respect to the Brownian sheet \( B(\cdot) \) over the region \( R_z \). Similarly,  \( \int_{R_z} \phi(\zeta) d\zeta \) represents the two-dimensional Lebesgue integral of the  function \( \phi \) over the region \( R_z \).  The process \( \{\phi(z), z \in R_z\} \) is a  measurable function of $ (\omega, z)$ and satisfies the condition $$ \int_{R_z} \mathbb{E}[\phi(\zeta)^2] \, d\zeta < \infty,$$ and for every \( \zeta \in R_z \), $ \phi(\zeta) $ is $ \mathcal{F}_{\zeta} $ -measurable. We denote the space of such processes by $\mathcal{L}^2_{a,1}$.
 
Let $\mathcal{L}^2_{a,2}$ denote the space of functions $ \phi(\omega, z, z')  $ defined on $ \Omega \times R_z \times R_z $. These functions are  measurable with respect $(\omega, \zeta, \zeta')$, and $\mathcal{F}_{\zeta \vee \zeta'} $-measurable for all $ \zeta, \zeta' \in R_z $. Additionally, the following condition is satisfied:
   \[
   \iint_{R_z \times R_z} I(\zeta \bar{\wedge} \zeta') E[\phi^2(\zeta,\zeta')] \, d\zeta \, d\zeta' < \infty.
   \]
Note that for two points $\zeta = (a_1, s_1)$ and $ \zeta' = (a_2, s_2) $, the operation $ \zeta \lor \zeta' $ is defined as the component-wise maximum of their coordinates: $$
\zeta \lor \zeta'  = (\max(a_1, a_2), \max(s_1, s_2)).$$
Further, for $a=(a_{1},a_{2})$ and $b=(b_{1},b_{2})$, we denote by $a%
\bar{\wedge}b$, if $a_{1}<b_{1}$ and $a_{2}>b_{2}$.  

It is important to note that the above double integrals are formulated so that only the values of the integrand at 
$\zeta \bar{\wedge} \zeta'$  contribute to each integral. Moreover, the indicator function $ I((t, x) \bar{\wedge} (t', x')) $ is defined as:
\[
I((t, x) \bar{\wedge} (t', x')) = 
\begin{cases} 
1, & \text{if } t \leq t' \text{ and } x \geq x', \\\\ 
0, & \text{otherwise}.
\end{cases}
\]
For any $ \phi \in \mathcal{L}^2_{a,2} $, the following stochastic integrals are well-defined (as shown in Wong \& Zakai (1978) \cite {WZ}) for all $ \zeta \in R_z$:
\[
\Phi(z) := \iint_{R_z \times R_z} \phi(\zeta,\zeta') \, B(d\zeta) \, B(d\zeta'),
\]
\[
\Phi_1(z) := \iint_{R_z \times R_z} \phi(\zeta,\zeta') \, d\zeta \, B(d\zeta'),
\quad
\Phi_2(z) := \iint_{R_z \times R_z} \phi(\zeta,\zeta') \, B(d\zeta) \, d\zeta'.
\]

In these cases $ \Phi$ is a martingale,  $\Phi_1 $ is an adapted 1-martingale,  
$ \Phi_2 $ is an adapted 2-martingale. Moreover, sample-continuous versions of these processes can be selected in all cases (see Wong \& Zakai (1978) \cite {WZ} for the definitions of adapted one- and two-martingales). 

Let \( f: \mathbb{R} \to \mathbb{R} \) be a continuously differentiable function, up to the fourth order (\( f \in C^4(\mathbb{R}) \)). The following theorem provides the It\^o formula for \( f(Y(z)) \), where \( Y(z) \) is a process driven by the Brownian sheet.

\begin{theorem}[It\^o Formula for \( f(Y(z)) \)]
Let \( Y(t,x) \) be a process as defined in (\ref{eq}), and $\alpha$ be continuous  function, $\beta \in \mathcal{L}^2_{a,1}  $, $\psi  \in \mathcal{L}^2_{a,2} $. Then the It\^o formula for $ f(Y(z))$ is given by:
\small
\begin{align*}\label{(17)}
&f(Y(z))  =f(Y_{0})+\int_{R_{z}}f^{\prime}(Y(\zeta))[\alpha(\zeta
)d\zeta+\beta(\zeta)B(d\zeta)]+\tfrac{1}{2}\int_{R_{z}}f^{\prime\prime}(Y(\zeta))\beta^{2}(\zeta
)d\zeta\nonumber\\
& +\iint\limits_{R_{z}\times R_{z}} \Big\{ f^{\prime\prime}(Y(\zeta\vee\zeta^{\prime
}))u \tilde{u}+f^{\prime} (Y(\zeta\vee \zeta^{\prime})) \psi(\zeta,\zeta^{\prime}) \Big\}B(d\zeta)B(d\zeta^{\prime}) \nonumber\\
&+\iint\limits_{R_{z}\times R_{z}}\Big\{f^{\prime\prime}(Y(\zeta\vee\zeta
^{\prime}))\Big( u \alpha(\zeta) + \psi(\zeta,\zeta^{\prime})  \tilde{u} \Big) +\tfrac{1}{2}f^{(3)}(Y(\zeta\vee\zeta^{\prime}))u^2 \tilde{u}\Big\}d\zeta B(d\zeta^{\prime})\nonumber\\
&+\iint\limits_{R_{z}\times R_{z}}\Big\{f^{\prime\prime}(Y(\zeta\vee\zeta
^{\prime}))\Big( \tilde{u}\alpha(\zeta^{\prime}) + \psi(\zeta,\zeta^{\prime}) u  \Big) +\tfrac{1}{2}f^{(3)}(Y(\zeta\vee\zeta^{\prime}))u^2 \tilde{u}\Big\}d\zeta B(d\zeta^{\prime})\nonumber\\
& +\iint\limits_{R_{z}\times R_{z}}I(\zeta \bar{\wedge} \zeta^{\prime}) \Big\{f^{\prime\prime}(Y(\zeta\vee\zeta
^{\prime}))\Big(\alpha(\zeta^{\prime})\alpha(\zeta)+\tfrac{1}{2}\psi^2(\zeta,\zeta^{\prime})\Big) \nonumber\\
&+f^{(3)}(Y(\zeta\vee\zeta^{\prime}))u \tilde{u}\psi(\zeta,\zeta^{\prime})+\tfrac{1}{2}f^{(3)}(Y(\zeta\vee\zeta^{\prime}))\left(\alpha(\zeta^{\prime}%
)\tilde{u}^2+\alpha(\zeta)u^2 \right)\nonumber\\
& + \tfrac{1}{4}f^{(4)}(Y(\zeta\vee\zeta^{\prime}))u^2 \tilde{u}^2\Big\}d\zeta d\zeta^{\prime},\nonumber
\end{align*}
where
$$u=\beta(\zeta^{\prime})+\int_{R_{z}}I(\zeta \bar{\wedge} \zeta^{\prime})\psi(\zeta,\zeta^{\prime})B(d\zeta),\,\,\,\textit{and} \,\,\,\,\tilde{u}=\beta(\zeta)+\int_{R_{z}}I(\zeta \bar{\wedge} \zeta^{\prime})\psi(\zeta,\zeta^{\prime})B(d\zeta^{\prime}).$$
$$$$
\end{theorem}
Note, that in the It\^o formula for \( f(Y(z)) \), the indicator function \( I(\zeta \bar{\wedge} \zeta') \) ensures proper symmetry and ordering within the double integral. Specifically, it enforces a consistent relationship between the variables \( t \) and \( x \) to maintain the correct integration structure.  

\subsection{Integration by Parts Formula}

Let \( Y_1(z) \) and \( Y_2(z) \) be two processes driven by Brownian sheets, defined as follows:
\[
Y_1(z) = Y_1(0) + \int_{R_{z}} \alpha_1(\zeta) \, d(\zeta) + \int_{R_{z}} \beta_1(\zeta) \, B(d(\zeta)) + \iint_{R_{z} \times R_{z}} \psi_1(\zeta, \zeta') \, B(d(\zeta)) B(d(\zeta')),
\]
\[
Y_2(z) = Y_2(0) + \int_{R_{z}} \alpha_2(\zeta) \, d(\zeta) + \int_{R_{z}} \beta_2(\zeta) \, B(d(\zeta)) + \iint_{R_{z} \times R_{z}} \psi_2(\zeta, \zeta') \, B(d(\zeta)) B(d(\zeta')).
\]
As shown  in Wong \& Zakai (1978) \cite{WZ}, the double  and the mixed integrals  $ \Phi$,  $\Phi_1 $, $ \Phi_2 $ are all weak martingales and therefore have expectation zero. Using  the It\^{o} formula above we get the following result.
\begin{lemma} \label{int_part}[Integration by Parts Formula]
The expectation \( \E[Y_1(z) Y_2(z)] \) is given by:
\small
\begin{align*}
  &  \E[Y_1(z)Y_2(z)]=Y_1(0)Y_2(0)\\
  &+\E\Big[\int_{R_{z}}\Big\{ Y_1(\zeta) \alpha_2(\zeta)+Y_2(\zeta) \alpha_1(\zeta)+ \beta_1(\zeta)\beta_2(\zeta)\Big\}d\zeta\nonumber\\
& +\iint\limits_{R_{z}\times R_{z}}I(\zeta \bar{\wedge} \zeta^{\prime}) \Big\{\alpha_1 (\zeta^{\prime}) \alpha_2(\zeta) + \alpha_{1}(\zeta)\alpha_2(\zeta^{\prime}) +\psi_1(\zeta,\zeta^{\prime}) \psi_2(\zeta,\zeta^{\prime}) \Big\} d\zeta d\zeta^{\prime}\Big].\\
\end{align*}
\end{lemma}
For further details and a proof, we refer to Agram et al. (2024) \cite{AOPT2}.

Let \( M = \{M(z), z \in R_{z_0}=[t_0, T] \times [x_0, X], T>0, X>0 \} \) be a square-integrable martingale adapted to the filtration \( \{\mathcal{F}_z\}_{z \in R_{z_0}} \). Then, there exists a unique pair of processes \( (\Psi, \Theta) \), where
\[
\Psi \in \mathcal{L}_{a,1}^2, \quad \Theta \in \mathcal{L}_{a,2}^2,
\]
such that
\[
M(z) = M(0) + \int_{R_z} \Psi(\zeta) dB(\zeta) + \int_{R_z \times R_z} \Theta(\zeta, \zeta') B(d\zeta) B(d\zeta'),
\]
for all \( z \in R_{z_0} \). Here:
\begin{itemize}
    \item The term $ \int_{R_z} \Psi(\zeta) dB(\zeta)$  represents the linear contribution of the Brownian sheet $ B $.
    \item The term \( \int_{R_z \times R_z} \Theta(\zeta, \zeta') B(d\zeta) B(d\zeta') \) captures second-order interactions between components of \( B \).
\end{itemize}

The following theorem from  Zaidi \& Nualart (2002)\cite{ZN} provides conditions for the existence and uniqueness of solutions to BSPDEs in the plane.
\begin{theorem} [Existence and Uniqueness]
Consider the BSPDE in the plane:
\begin{align*}
Y_z &= \xi - \int_{R_{z_0} \setminus R_z} f(\zeta, Y(\zeta), \Psi(\zeta)) d\zeta - \int_{R_{z_0} \setminus R_z}  \Psi(\zeta) dB(\zeta) \\&- \int_{R_{z_0} \setminus R_z}\int_{R_{z_0} \setminus R_z}                  \Theta(\zeta, \zeta') B(d\zeta) B(d\zeta').
\end{align*}
Assume the following conditions:
\begin{enumerate}
    \item \( f(z, y, \Psi) \) is Lipschitz in \( y \) and \( \Psi \), with constants \( K_1 \) and \( K_2 \) with Lipschitz constants satisfy:
    \[
    K_1 |z_0| < \sqrt{r_0}, \quad K_2^2 |z_0| < 1,
    \]
    where \( r_0 \) is the first positive root of the Bessel function \( J_0(2\sqrt{t}) \).
    \item The terminal condition \( \xi \) satisfies:
    \[
    \xi \in L^2(\Omega, \mathcal{F}_{z_0}, \mathbb{P}).
    \]
\end{enumerate}
There exists a unique triple of processes \( (Y, \Psi, \Theta) \), with:
\[
Y \in \mathcal{L}^2_{a,1}, \quad \Psi \in \mathcal{L}^2_{a,1}, \quad \Theta \in \mathcal{L}^2_{a,2},
\]
such that \( Y_z \) satisfies the BSPDE for all \( z \leq z_0 \).
\end{theorem}


\section{Nonzero-sum stochastic differential games}
This section focuses on a stochastic differential game involving two players in  a nonzero-sum setting, where the dynamics of the state process are  governed a stochastic reaction-diffusion equation   with spatial and temporal dependencies. The problem is formulated in the domain $\R \times \R^n$ and the state evolves according to the following equation:
\begin{align}
\frac{\partial Y(t,x)}{\partial t} = F(t,x,Y(t,x),u(t,x)) + \sigma(t,x) \frac{\partial^2 B(t,x)}{\partial t \partial x}.
\end{align}
$F$ is a  function representing system dynamics, $\sigma(t,x) \frac{\partial^2 B(t,x)}{\partial t \partial x}$ introduces spatial-temporal stochastic noise.

The controls $u_i(t,x)$, $i=1,2$ represent the strategies of the two players in a nonzero-sum stochastic differential game. Each player aims to optimize the individual performance functional, which depends on the state process $Y(t,x)$, the controls, and the time-space variables $(t,x)$.

Let $U_i\subset\mathbb{R}$, $i=1,2$ be  the set of  allowable control values for player $i$, and let $\mathcal{U}_i$ denote the set of all $\mathcal{F}_{t,x}$-adapted control processes $u_i(t,x)$, $t<T$, $x<X$ valued in $U_i$.  We therefore define the set of admissible control processes $\mathcal{A}_i\subset\mathcal{U}_i$ to be the collection of all $\mathcal{F}_{t,x}$-adapted processes with values in $U_i$.
The controls $u_1(t,x)$ and $u_2(t,x)$ determine the evolution of the state process $Y(t,x)$, and each player's strategy impacts the overall system dynamics.


The state process $Y(t,x)$ evolves according to some SPDE with spatial and temporal dependencies. Specifically, the state $Y(t, x)$ is driven by the controls of the players $u_i(t,x)$, $i=1,2$  and  some noise. Each player's control impacts the system, and their objective is to maximize their respective performance functionals. Let us denote $u=(u_1,u_2)$ and similarly $U=(U_1,U_2)$ and  $\mathcal{U}=(\mathcal{U}_1,\mathcal{U}_2)$.

Let the set of admissible control processes for Player I and Player II are denoted by  $\mathcal{A}_1$ and $\mathcal{A}_2$, respectively. The game is formulated in continuous time, and the cost functionals for Player I and Player II are defined as follows:
\begin{itemize}
    \item {Player I's performance functional}:
    \[
    J_1(u_1, u_2) = \mathbb{E} \left[ \int_{{R}_Z} \left( f_1(\zeta, Y(\zeta), u_1(\zeta), u_2(\zeta)) \right) d\zeta + g_1(Y(Z)) \right],
    \]
    \item {Player II's performance functional}:
    \[
    J_2(u_1, u_2) = \mathbb{E} \left[ \int_{{R}_Z} \left( f_2(\zeta, Y(\zeta), u_1(\zeta), u_2(\zeta)) \right) d\zeta + g_2(Y(Z)) \right],
    \]
\end{itemize}
where $R_Z=[0,T] \times [0,X],$ represents the time-space interval over which the game takes place, and $Z=(T,X)$ is the terminal time-space, for  given $T>0,\, X>0$.
For  $i=1,2$, the function $f_i:[0, T] \times [0, X]\times \mathbb{R} \times U \to \mathbb{R}$ is an $\mathcal{F}_{t,x}$-adapted and $C^1$ with respect to $y$ and $u_1$. It represents the running cost for player $i$, which depends on the state $Y$, controls $u_1$ and $u_2$, and the time-space variable $(t,x)$, $g_i:[0, T] \mathbb{R} \to \mathbb{R}$ represents the terminal cost which is an $\mathcal{F}_{T,X}$-measurable square integrable random variable and $C^1$ with respect to $y$, depending on the state at terminal time $Z$. The functions $f_1$, $f_2$, $g_1$, and $g_2$ are generally convex in $u_1$ and $u_2$, ensuring that each player's optimization problem is well-posed.

The dynamics of the state $Y(z)=Y^u(z)$ for  $z\leq Z$ of the system is described by the equation
\small
\begin{align}\label{H_eq}
Y(z)=Y(t,x)=Y(0)+\int_{R_{z}}\alpha(\zeta,Y(\zeta),u(\zeta))d\zeta+\int_{R_{z}}\beta(\zeta,Y(\zeta),u(\zeta))B(d\zeta),
\end{align}
where $R_z=[0,t] \times [0,x]$ when $z=(t,x)$, and $u$ denotes a control process. The coefficients $\alpha(t,x,Y,u), \beta(t,x,Y,u):[0, T] \times [0, X]\times \mathbb{R} \times U \to \mathbb{R}$ are $\mathcal{F}_{t,x}$-adapted and $C^1$ with respect to $y$ and $u$. 

The intuitive idea is that there are two players, I and II,  each of the players is trying to maximize their respective performance functionals.\\
A pair of controls $(\widehat{u}_1,\widehat{u}_2) \in \A_1 \times \A_2$ is called a \emph{Nash equilibrium} if it satisfies the following conditions:
\begin{equation}\label{P1}
    J_1(\widehat{u}_1,\widehat{u}_2) \leq J_1(u_1, \widehat{u}_2) \text{ for all  } u_1 \in \A_1
\end{equation}
and 
  \begin{equation}\label{P2}
    J_2(\widehat{u}_1,\widehat{u}_2) \leq J_2(\widehat{u}_1,u_2) \text{ for all  } u_2 \in \A_2.
\end{equation}
\begin{figure}[!ht]
\center{\includegraphics[scale=0.37]{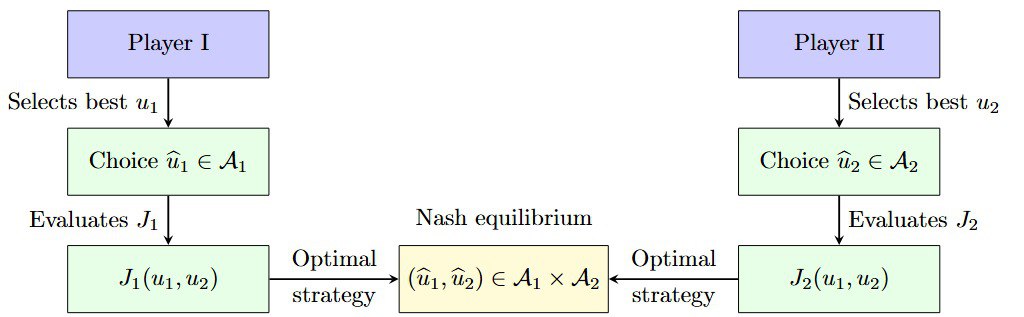}}
\caption{Decision process leading to Nash equilibrium in a two-player game model.}
\label{fig}
\end{figure}

In other words, as long as Player II using the control $\widehat{u}_2$ the best for Player I is to use $\widehat{u}_1$ and vice versa.
Therefore a Nash equilibrium, although it is not in general unique, is in a sense a likely outcome of a game. It can in some cases be regarded as a fixed point of the game and hence a limit strategy pair if the game is played iteratively.

The optimality principle provides first-order necessary conditions for optimality in this stochastic differential game  framework. For a pair of controls $(\widehat{u}_1,\widehat{u}_2)$  to be optimal, there must exist adjoint processes satisfying backward stochastic partial differential equations (BSPDEs) driven by Brownian sheet, ensuring that the Hamiltonian of the system is minimized pointwise with respect to the controls. This principle serves as a critical tool for analyzing and characterizing the optimal strategies for both players. Our result may be regarded as an extension of the results in Agram et al. (2025) \cite{AOPT1}. Furthermore, SPDEs of the type \eqref{H_eq} have been extensively studied over the years by various authors. Notably, Cairoli (1972)\cite{Cairoli72} and Yeh (1981)\cite{Yeh81}  established strong existence and pathwise uniqueness of solutions to equations of the form \eqref{H_eq}. For the case of weak solutions, we also refer to Yeh (1985)\cite{Yeh85} .

Suppose there exists a Nash equilibrium $(\widehat{u}_1,\widehat{u}_2)$.
To find it, based on the optimality principle for stochastic control. To do this we introduce the following:
\begin{definition}[Convolution-like operator]\label{Clo}
For general functions $h,k:[0,T]\times[0,X] \to \mathbb{R}$, we define
\begin{align*}
 (h\star k)(t,x)=\int_0^x {\int_{t}^T }\{ h(t,x)k(s,a) + h(s,a) k(t,x)\} ds da.  
\end{align*}
\end{definition}
Note that with this notation, we have
\small
\begin{align*}
 \iint\limits_{R_{z}\times R_{z}}I(\zeta \bar{\wedge} \zeta^{\prime}) \Big\{\alpha_1 (\zeta^{\prime}) \alpha_2(\zeta) + \alpha_{1}(\zeta)\alpha_2(\zeta^{\prime})  \Big\} d\zeta d\zeta^{\prime}=\int_{R_{z}}( \alpha_1 \star \alpha_2)(\zeta) d\zeta.  
\end{align*}

For each $i=1,2$, define the associated Hamiltonian 
$$H_i: R_{Z} \times \R \times U_1 \times U_2\times \R \times \R \times \R^{R_{Z}}  \to \R$$ as follows:
\begin{align}
 & H_i(z, y, u_1, u_2 ,p_i,q_i, L_i)\nonumber\\
 & :=f_i(z,y,u_1, u_2)+\alpha(z,y,u_1, u_2)p_i+\beta(z,y,u_1, u_2)q_i+(L_i \star \alpha)(z).
\end{align}
Here $\R^{R_{Z}}$ is a set of functions from $R_{Z}$ to $\R$. \\ Following, Definition \ref{Clo}, the convolution-like operator $(L_i \star \alpha)(z)$ is defined by
$$ 
(L_i \star \alpha)(z)=\int_0^x\int_t^T \Big\{L_i(z) \alpha(\zeta',y,u_1, u_2) + L_i(\zeta')\alpha(z,y,u_1, u_2)\Big\}d\zeta' ;\quad z=(t,x),$$ 
where $L_i =L_i(z) $ and the adjoint processes $(p_i,q_i, r_i)=(p_i(z),q_i(z), r_i(z, \cdot)$ are given by the following equations:\\
$L_i:R_Z \mapsto \mathbb{R}$ is defined implicitly by the integral equation
\begin{align}
  L_i(z)= -\frac{\partial H_i}{\partial y} (z,Y(z),u(z), p_i(z),q_i(z),r_i(z,\cdot),L_i(z)); \quad z=(t,x)
\end{align}

and $p_i,q_i,r_i$ are defined as the solution of the BSPDE
\begin{align}
\begin{cases}
    p_i(dz)&=- \frac{\partial H_i}{\partial y}(z)dz+q_i(z) B(dz) 
     +M_i^{r_i}(dz),\quad 0 \leq (t,x) \leq(T,X),\\
    p_i(Z)&=\frac{\partial g_i}{\partial y}(Y(Z)).
    \end{cases}
    \end{align}
   or, in integrated form,
   \small
    \begin{align}
   p_i(z)&=\frac{\partial g_i}{\partial y}(Y(Z))- \int_{R_{Z} \setminus R_z}\frac{\partial H_i}{\partial y}(\zeta
   ) d\zeta+\int_{R_{Z} \setminus R_z}\ q_i(\zeta)B(d\zeta)
   + M_{i}^{r_i}(Z)- M_{i}^{r_i}(z),\quad z\leq Z, 
    \end{align}
where the term  $M_{i}^{r_i}, i=1,2$ represents a second-order stochastic integral 
\begin{align*}
    M_{i}^{r_i}(z)=\iint\limits_{R_{Z}\times R_{Z}} r_i(\zeta,\zeta') B(d\zeta)B(d\zeta'),\quad r_i \in \mathcal{L}_{a,2}, \quad i=1,2.
\end{align*}

\begin{remark} 
(i) Since \( H_i \) itself depends on \( L_i \), the equation
\begin{align}
L_i(z) = -\frac{\partial H_i}{\partial y} (z,Y(z),u(z), p_i(z),q_i(z),r_i(z,\cdot),L_i(z))
\end{align}
forms a self-consistent integral equation. It follows that \( L_i \) satisfies an integral equation of Volterra type.\\
(ii) The adjoint variables $L_i$, obtained as a solution to an integral equation, are new. They come in addition to the adjoint variables $(p_i,q_i,r_i)$ in the corresponding maximum principle for the classical (1-parameter) Brownian motion.\\
(iii) Note that the operators $ L_i \star \alpha$ are bilinear; \ i=1,2. 

\end{remark}

There are two versions of the optimality principle for this problem, namely the so-called \emph{sufficient} and \emph{necessary} conditions of optimality.
We present them both below.

\subsection{Sufficient Optimality Conditions for Two Players}
 We formulate the sufficient optimality condition, which is particularly useful for verifying optimality, since it specifies the conditions under which no alternative strategy leads to lower costs for any of the players.

\begin{theorem}[Sufficient optimality condition]
Let $(\widehat{u}_1(\zeta), \widehat{u}_2(\zeta)) \in \A_1 \times \A_2$ be admissible controls for Player I and Player II, respectively, with corresponding state process $\widehat{Y}(\zeta)$, and adjoint processes $(\widehat{p}_i(\zeta), \widehat{q}_i(\zeta), \widehat{r}_i(\zeta), \widehat{L}_i(\zeta))$ for each  player $i=1,2$. Suppose these processes satisfy the following conditions:
\begin{enumerate}
    \item The functions $f_1(\zeta, Y, u_1, u_2)$ and $f_2(\zeta, Y, u_1, u_2)$ are convex in both $u_1$ and $u_2$, and the terminal cost functions $g_1(Y)$ and $g_2(Y)$ are convex in $Y$.
    \item The Hamiltonians $H_1(\zeta, Y, u_1, u_2, p_1, q_1, L_1)$ and $H_2(\zeta, Y, u_1, u_2, p_2, q_2, L_2)$ are convex in both $u_1$ and $u_2$, for all $(Y, p_1, q_1, L_1, p_2, q_2, L_2)$.
    \item The following conditions hold for the control pair $(\widehat{u}_1(\zeta), \widehat{u}_2(\zeta))$:
    \begin{itemize}
        \item For Player I:
    \end{itemize}
       \begin{align*}
            \min_{v_1 \in \mathcal{A}_1} H_1(\zeta, \widehat{Y}(\zeta), v_1, \widehat{u}_2(\zeta), \widehat{p}_1(\zeta), &\widehat{q}_1(\zeta), \widehat{L}_1(\zeta))&\\ = &H_1(\zeta, \widehat{Y}(\zeta), \widehat{u}_1(\zeta), \widehat{u}_2(\zeta), \widehat{p}_1(\zeta), \widehat{q}_1(\zeta), \widehat{L}_1(\zeta)),
       \end{align*}
            \begin{itemize}
            \item  For Player II:
        \end{itemize}
             \begin{align*}
        \min_{v_2 \in \mathcal{A}_2} H_2(\zeta, \widehat{Y}(\zeta), \widehat{u}_1(\zeta), v_2, \widehat{p}_2(\zeta), &\widehat{q}_2(\zeta), \widehat{L}_2(\zeta))&\\ = &H_2(\zeta, \widehat{Y}(\zeta), \widehat{u}_1(\zeta), \widehat{u}_2(\zeta), \widehat{p}_2(\zeta), \widehat{q}_2(\zeta), \widehat{L}_2(\zeta)).
        \end{align*}
   
\end{enumerate}
Then $(\widehat{u}_1(\zeta), \widehat{u}_2(\zeta))$ are optimal strategies for both players.
\end{theorem}

\dproof
Let $(\widehat{u}_1, \widehat{u}_2)$ satisfy the conditions above.  The equilibrium conditions require that each player cannot improve his position by changing his strategy given the other player's strategy. It is important that the inequalities are non-strict, meaning that the cost associated with a player's current strategy cannot be lower than the cost of any alternative strategy available to another player.
We now proceed to show that the controls $ (\widehat{u}_1, \widehat{u}_2)$ are optimal by proving that any deviation from these controls results in a higher cost for at least one player.

The difference in the cost functionals when deviating from $\widehat{u}_1$   for Player I can be written as:
\[
J_1(u_1,  \widehat{u}_2) - J_1(\widehat{u}_1, \widehat{u}_2) = I_1 + I_2,
\]
where
\[
I_1 = \mathbb{E} \left[ \int_{{R}_Z} \left( f_1(\zeta, Y(\zeta), u_1(\zeta),  \widehat{u}_2) - f_1(\zeta, \widehat{Y}(\zeta), \widehat{u}_1(\zeta), \widehat{u}_2(\zeta)) \right) d\zeta \right],
\]
\[
I_2 = \mathbb{E} \left[ g_1(Y(Z)) - g_1(\widehat{Y}(Z)) \right].
\]
The term $I_1$  represents the accumulated difference in running costs, while  $I_2$  takes into account the difference in final costs.
Now express  $I_1$  using the definition of the Hamiltonian function $H_1$, which represents the dynamics of the system and the costs associated with the players' controls:
\begin{align}\label{I_1} 
I_1 &= \mathbb{E} \int_{{R}_Z}\Big[  H_1(\zeta, Y(\zeta), u_1(\zeta), \widehat{u}_2(\zeta), p_1(\zeta), q_1(\zeta), L_1(\zeta)) \nonumber\\&- H_1(\zeta, \widehat{Y}(\zeta), \widehat{u}_1(\zeta), \widehat{u}_2(\zeta),\widehat{p}_1(\zeta), \widehat{q}_1(\zeta), \widehat{L}_1(\zeta)) \nonumber 
\\&-\widehat{p}_1(\zeta) \tilde{\alpha}(\zeta)-\widehat{q}_1(\zeta)\tilde{\beta}(\zeta) -(\widehat{L}_1 \star \tilde{\alpha})(\zeta) \Big] d\zeta,
\end{align}
where $\tilde{\alpha}(\zeta)=\alpha(\zeta,{Y}(\zeta), u_1(\zeta),\widehat{u}_2(\zeta))- \alpha(\zeta,\widehat{Y}(\zeta),\widehat{u}_1(\zeta),\widehat{u}_2(\zeta))$ etc.
For $I_2$, the convexity of $g_1$ gives
$$g_1(Y(Z)) - g_1(\widehat{Y}(Z))\geq \frac{\partial g_1}{\partial y} (\widehat{Y}(Z))(Y-\widehat{Y}(Z)),$$
or
$$g_1(Y(Z)) - g_1(\widehat{Y}(Z))\geq \frac{\partial g_1}{\partial y} (\widehat{Y}(Z))\tilde{Y}(Z).$$
Thus, using the integration by part Lemma \ref{int_part}, the fact that the $B(dz)$-integrals and the $B(dz)B(dz')$-integrals are orthogonal (see Cairoli \& Walsh (1975) \cite{CW}, Theorem 2.5), implies
\begin{align}
    I_2 &\geq \E\Big[\frac{\partial g_1}{\partial y} (\widehat{Y}(Z))\tilde{Y}(Z)\Big]=\E\Big[\widehat{p}_1(Z)\tilde{Y}(Z)\Big]\nonumber\\
    &=\E\Big[\int_{R_{Z}}\Big\{\widehat{p}_1(\zeta)\tilde{\alpha}(\zeta)
    - \frac{\partial \hat{H}_1}{\partial y}\tilde{Y} (\zeta)+\widehat{q}_1(\zeta)\tilde{\beta}(\zeta) \textcolor{blue}{+}(\widehat{L}_1\star  \tilde{\alpha})(\zeta)\Big\}d\zeta\Big].\label{4.8}
\end{align} 
Adding \eqref{I_1}, \eqref{4.8} and using the assumption that the Hamiltonian is convex in $u_1$,  we obtain 
\begin{align*}
   J_1(u_1,  \widehat{u}_2) - J_1(\widehat{u}_1, \widehat{u}_2) &=I_1 + I_2 \geq \E\Big[\int_{R_{Z}}\Big\{ H_1(\zeta)-\widehat{H}_1(\zeta) -\frac{\partial \hat{H}_1}{\partial y}(\zeta)\tilde{Y}(\zeta)\Big\}d\zeta\Big] \nonumber\\ 
    &\geq \E\Big[\int_{R_{Z}} \frac{\partial \hat{H}_1}{\partial u_1} (\zeta)\tilde{u}_1(\zeta)d\zeta\Big] \geq 0.
\end{align*}
This proves that 
\[
J_1(u_1,  \widehat{u}_2) - J_1(\widehat{u}_1, \widehat{u}_2)=I_1+I_2 \geq 0.
\]
and therefore $\widehat{u}_1$ is optimal. This shows that deviating from $\widehat{u}_1$ cannot decrease Player I's cost.

By symmetry and similar reasoning for Player II, we also have:
\[
J_2(\widehat{u}_1,  u_2) - J_2(\widehat{u}_1, \widehat{u}_2) \geq 0.
\]

Since any deviation from the optimal control pair $ (\widehat{u}_1, \widehat{u}_2)$ does not decrease the cost for either player, we conclude that $ (\widehat{u}_1, \widehat{u}_2)$
is indeed the optimal control pair, corresponding to the definition of a Nash equilibrium.
\fproof

\begin{remark}  
The inequality  
$ J_1 (u_1, \widehat{u}_2) - J_1 (\widehat{u}_1, \widehat{u}_2) \geq 0$  
is one Nash equilibrium condition given in \eqref{P1}, meaning that Player I cannot reduce their cost by deviating from \( \widehat{u}_1 \), given that Player II follows \( \widehat{u}_2 \).  
Similarly, for Player II, the inequality  
$
J_2(\widehat{u}_1,  u_2) - J_2(\widehat{u}_1, \widehat{u}_2) \geq 0$  
is the other equilibrium condition \eqref{P2}, indicating that Player II also cannot lower their cost by unilaterally changing their strategy from \( \widehat{u}_2 \) while Player I follows \( \widehat{u}_1 \).  
\end{remark}



\subsection{Necessary Optimality Conditions for Two Players}

In the context of a two-player non-zero sum game, the problem setup involves two players, each optimizing their own cost functional while interacting with the other player through their control variables. The goal is to extend the necessary condition of optimality to account for interactions between the two players and their respective controls. The framework we consider involves two control processes, \( u_1 \) and \( u_2 \), corresponding to Player~I and Player~II, respectively, and their dynamics are influenced by both their own control and the opponent's control.

Let \( u_1 \) and \( u_2 \) be the controls of Player I and Player II, respectively, and the state dynamics depend on both controls. The cost functionals for the two players are given by:
\[
J_1(u_1, u_2) = \mathbb{E}\left[ \int_{R_Z} f_1(\zeta, Y^{(u_1, u_2)}(\zeta), u_1(\zeta),  u_2(\zeta), p(\zeta), q(\zeta)) \, d\zeta + g_1(Y^{(u_1, u_2)}(Z)) \right]
\]
and
\[
J_2(u_1, u_2) = \mathbb{E}\left[ \int_{R_Z} f_2(\zeta, Y^{(u_1, u_2)}(\zeta), u_1(\zeta), u_2(\zeta), p(\zeta), q(\zeta) ) \, d\zeta + g_2(Y^{(u_1, u_2)}(Z)) \right],
\]
where $Y^{(u_1, u_2)}$ denotes the state variable, evolving according to a stochastic differential equation influenced by both players' controls, the functions \( f_1 \), \( f_2 \), \( g_1 \), and \( g_2 \) represent the running cost and terminal cost for Players I and II, $p$ and $q$ denote additional state-dependent terms that contribute to the evolution of the dynamics and the cost functions.

We assume the following conditions on the set \( \mathcal{A}=\A_1 \times \A_2 \) of admissible controls for both players:

\begin{itemize}
    \item \textbf{(A1)}: The sets  \( \mathcal{A}_1 \) and \( \mathcal{A}_2 \) of admissible controls are convex. This ensures that linear combinations of admissible controls for both players remain admissible.
    \item \textbf{(A2)}: For any \( z_0 = (t_0, x_0) < Z = (T, X) \) and any bounded \( \mathcal{F}_{z_0} \)-measurable random variable \( \eta_{z_0} \), the control pairs $u_{1, z_0}(\zeta)$ and $u_{2, z_0}(\zeta)$  are admissible for Player I and Player II, where
    \[
    u_{1, z_0}(\zeta) = \eta_{z_0} \mathbf{1}_{R_{z_0}}(\zeta) \quad \text{and} \quad u_{2, z_0}(\zeta) = \eta_{z_0} \mathbf{1}_{R_{z_0}}(\zeta)
    \]
  with \( \mathbf{1}_{R_{z_0}}(\zeta) \) being the indicator function of the rectangle \( R_{z_0} = [t_0, T] \times [x_0, X] \). This condition allows for specific perturbations in the controls of both players.
\end{itemize}

We consider the effect of small perturbations in both players' controls and derive the necessary conditions.

\subsection{Derivative Process for Two Players}
To derive the necessary optimality conditions, we consider small perturbations in the controls of both players. Specifically, for a small positive value of $\epsilon$ that approaches zero, we introduce small perturbations $\epsilon v_1$ in Player I's control $u_1$, keeping Player II's control $u_2$ fixed. Similarly, we introduce a small perturbation $\epsilon v_2$ in Player II's control $u_2$, keeping Player I's control $u_1$ fixed.

To examine the impact of perturbations in the controls of both players, we define the derivative processes \( G_1(\zeta) \) and \( G_2(\zeta) \), which capture the change in state variables due to small perturbations in Player's I control $u_1$ and Player's II control $u_2$ 
\[
G_1(\zeta) := \lim_{\epsilon \to 0} \frac{1}{\epsilon} \left( Y^{(u_1 + \epsilon v_1, u_2)}(\zeta) - Y^{(u_1, u_2)}(\zeta) \right)
\]
and
\[
G_2(\zeta) := \lim_{\epsilon \to 0} \frac{1}{\epsilon} \left( Y^{(u_1, u_2 + \epsilon v_2)}(\zeta) - Y^{(u_1, u_2)}(\zeta) \right).
\]

These derivative processes, also known as G\^ateaux derivatives (see  \O ksendal (2003)\cite{sde} or Fleming \& Rishel (1975)\cite{Fl})  quantify the change in the state variables of each player due to small perturbations in the control of the other player.

The equations for \( G_1(\zeta) \) and \( G_2(\zeta) \) satisfy:
\begin{align*}
&G_1(z)-G_1(0) \\&= \int_{R_{z_0}} \left\{ \frac{\partial \alpha_1}{\partial y}(\zeta) G_1(\zeta) + \frac{\partial \alpha_1}{\partial u_1}(\zeta) v_1(\zeta) \right\} d\zeta
+ \int_{R_{z_0}} \left\{ \frac{\partial \beta_1}{\partial y}(\zeta) G_1(\zeta) + \frac{\partial \beta_1}{\partial u_1}(\zeta) v_1(\zeta) \right\} B(d\zeta),
\end{align*}
and
\begin{align*}
&G_2(z)-G_2(0) \\&
 =  \int_{R_{z_0}} \left\{ \frac{\partial \alpha_2}{\partial y}(\zeta) G_2(\zeta) + \frac{\partial \alpha_2}{\partial u_2}(\zeta) v_2(\zeta) \right\} d\zeta
+ \int_{R_{z_0}} \left\{ \frac{\partial \beta_2}{\partial y}(\zeta) G_2(\zeta) + \frac{\partial \beta_2}{\partial u_2}(\zeta) v_2(\zeta) \right\} B(d\zeta).
\end{align*}
These equations describe the evolution of the perturbation processes for both players, which depend on their own control and the control of the opponent. The detailed formulation and analysis of these processes $G_i, i=1,2$, including the associated equations, are provided in Agram et al. (2025) \cite{AOPT1}

The variation of the cost functionals for both players due to small perturbations in their controls is given by:
\[
\frac{d}{d\epsilon} J_1(u_1 + \epsilon v_1, u_2) \Big|_{\epsilon = 0} = \mathbb{E}\left[ \int_{R_Z} \frac{\partial H_1}{\partial u_1}(\zeta, Y^{(u_1, u_2)}(\zeta), u_1(\zeta), u_2(\zeta), p_1(\zeta), q_1(\zeta)) v_1(\zeta) \, d\zeta \right]
\]
and
\[
\frac{d}{d\epsilon} J_2(u_1, u_2 + \epsilon v_2) \Big|_{\epsilon = 0} = \mathbb{E}\left[ \int_{R_Z} \frac{\partial H_2}{\partial u_2}(\zeta, Y^{(u_1, u_2)}(\zeta), u_1(\zeta), u_2(\zeta), p_2(\zeta), q_2(\zeta)) v_2(\zeta) \, d\zeta \right].
\]

The basic idea behind the necessary optimality condition is that for each player, the variation of his cost functional must be non-positive under all admissible perturbations of the opponent's control. This ensures that no player can improve his costs by unilaterally changing the control, provided that the other player's strategy remains unchanged.
\begin{theorem}
\begin{itemize}
    \item \textbf{Necessary Condition for Player I}: Suppose \( \widehat{u}_1 \in \mathcal{A}_1 \) is the optimal control for Player I, and \( \widehat{u}_2 \in \mathcal{A}_2 \) is the optimal control for Player II. Then the necessary condition for optimality for Player I is:
    \[
    \frac{\partial H_1}{\partial u_1} (\zeta, \widehat{Y}^{(u_1, u_2)}(\zeta), \widehat{u}_1(\zeta), \widehat{u}_2(\zeta), \widehat{p}_1(\zeta), \widehat{q}_1(\zeta), \widehat{L}(\zeta)) = 0 \quad \text{for almost all} \quad \zeta.
    \]
    \item \textbf{Necessary Condition for Player II}: Similarly, the necessary condition for optimality for player II is:
    \[
    \frac{\partial H_2}{\partial u_2} (\zeta, \widehat{Y}^{(u_1, u_2)}(\zeta), \widehat{u}_1(\zeta),\widehat{u}_2(\zeta), \widehat{p}_2(\zeta), \widehat{q}_2(\zeta), \widehat{L}(\zeta)) = 0 \quad \text{for almost all} \quad \zeta.
    \]
\end{itemize}
\end{theorem}
Thus, the necessary optimality condition for a two-player game provides the conditions under which both players' controls are optimal in response to each other's strategies. Each player's control must satisfy the condition that the derivative of the Hamiltonian with respect to their control vanishes almost everywhere, ensuring that the strategy of each player is optimal given the strategy of the other player.
\newline

\dproof Consider the variation of the functional $J_1$ with respect to the control $u_1$.
We compute the derivative:
\begin{equation*}
    \frac{d}{d\epsilon} J_1(u_1 + \epsilon v_1, u_2) \Big|_{\epsilon=0}.
\end{equation*}
Expanding the objective functional:
\begin{align*}
    &\frac{d}{d\epsilon} J_1(u_1 + \epsilon v_1, u_2) \Big|_{\epsilon=0} =\lim_{\epsilon \rightarrow 0} \frac{1}{\epsilon} \mathbb{E} \Bigg[\int_{R_{Z}} \Big\{ f_1(\zeta, Y^{(u_1+\epsilon v_1, u_2)}(\zeta), u_1(\zeta)+\epsilon v_1(\zeta), u_2(\zeta)) \\
    &\quad - f_1(\zeta, Y^{(u_1, u_2)}(\zeta), u_1(\zeta), u_2(\zeta)) \Big\} d\zeta + g_1(Y^{(u_1+\epsilon v_1, u_2)}(Z)) - g_1(Y^{(u_1, u_2)}(Z)) \Bigg].
\end{align*}
Using the first-order expansion:
\begin{align*}
   \frac{d}{d\epsilon} J_1(u_1 + \epsilon v_1, u_2) \Big|_{\epsilon=0}  &= \mathbb{E} \Bigg[\int_{R_{Z}} \Big\{\frac{\partial f_1}{\partial y} (\zeta, Y^{(u_1, u_2)}(\zeta), u_1(\zeta), u_2(\zeta)) G_1(\zeta) \\
    &\quad +\frac{\partial f_1}{\partial u_1} (\zeta, Y^{(u_1, u_2)}(\zeta), u_1(\zeta), u_2(\zeta)) v_1(\zeta)\Big\} d\zeta  + \frac{\partial g_1}{\partial y} (Y^{(u_1, u_2)}(Z)) G_1(Z) \Bigg], 
\end{align*}
where $G_1$ is the variation in the system's state due to the change in control.
We now decompose it into two terms:
\[
 \frac{d}{d\epsilon} J_1(u_1 + \epsilon v_1, u_2) \Big|_{\epsilon=0}=I_1 + I_2,
\]
where
\begin{align*}
  I_1 &= \mathbb{E} \Bigg[\int_{R_{Z}} \Big( \frac{\partial H_1}{\partial y}(\zeta) -\frac{\partial \alpha}{\partial y} (\zeta)p_1(\zeta)-\frac{\partial \beta}{\partial y} (\zeta) q_1(\zeta)-(L_1 \star  \frac{\partial \alpha}{\partial y})(\zeta)\Big) G_1(\zeta) d\zeta \\
     &\quad +\int_{R_{Z}}  \Big( \frac{\partial H_1}{\partial u_1} (\zeta) -\frac{\partial \alpha}{\partial u_1} (\zeta)p_1(\zeta)-\frac{\partial \beta}{\partial u_1} (\zeta)  q_1(\zeta)-(L_1 \star \frac{\partial \alpha}{\partial u_1})(\zeta)\Big) v_1(\zeta)  d\zeta\Bigg],
\end{align*}
and
\begin{align*}
I_2 &= \mathbb{E} \Bigg[\frac{\partial g_1}{\partial y} (Y^{(u_1, u_2)}(Z)) G_2(Z) \Bigg] = \mathbb{E} \Big[p_1 (Z) G_2(Z)\Big]\\
&= \mathbb{E} \Bigg[\int_{R_{Z}}\Big\{ p_1(\zeta)\Big\{ \frac{\partial \alpha}{\partial y} (\zeta) G_2(Z) + \frac{\partial \alpha}{\partial u_1}(\zeta)v_1(\zeta)\Big\} - \frac{\partial H_1}{\partial y} (\zeta)G_2(Z)\\
&\quad + q_1(\zeta) \Big\{ \frac{\partial \beta}{\partial y}(\zeta) G_2(\zeta)+\frac{\partial \beta}{\partial u_1}(\zeta)v_1(\zeta)\Big\}-\big\{ \frac{\partial H_1}{\partial y} \star (\frac{\partial \alpha}{\partial y} G_2+ \frac{\partial \alpha}{\partial u_1}  v_1)\big\} (\zeta)\Big\}d\zeta\Bigg].
\end{align*}
Adding \( I_1 \) and \( I_2 \), and choosing \( L_1 = - \frac{\partial H_1}{\partial y} \), we obtain:
\[
\frac{d}{d\epsilon} J_1(u_1+\epsilon v_1, u_2) \Big|_{\epsilon=0} = \mathbb{E} \Bigg[\int_{R_{Z}} \frac{\partial H_1}{\partial u_1}(\zeta) v_1(\zeta)  d\zeta\Bigg].
\]
For a test function of the form:
\[
v_1(\zeta)=\eta_{z_{0}} \mathbf{1}_{R_{z_{0}}}(\zeta),
\]
this leads to
\[
\mathbb{E} \Bigg[\int_{R_{z_0}} \frac{\partial H_1}{\partial u_1}(\zeta)\eta_{z_{0}} d\zeta\Bigg]\leq 0.
\]
Since this holds for all \( z_0 \), we deduce:
\[
\frac{\partial ^2}{\partial t_0 \partial x_0}  \left( \mathbb{E} \left[\int_{R_{z_0}} \frac{\partial H_1}{\partial u_1}\eta_{z_{0}} d\zeta\right]\right)=\frac{\partial H_1}{\partial u_1}(z_{0}) \eta_{z_{0}} \leq 0.
\]
Since this holds for all \( \eta_{z_{0}} \in\mathbb{R} \), we conclude:
\[
\frac{\partial H_1}{\partial u_1}(z_{0})=0.
\]
A similar derivation follows for Player II. By performing the same expansion for \( J_2 \), we obtain:
\[
\frac{d}{d\epsilon} J_2(u_1, u_2+\epsilon v_2) \Big|_{\epsilon=0} = \mathbb{E} \Bigg[\int_{R_{Z}} \frac{\partial H_2}{\partial u_2}(\zeta) v_2(\zeta)  d\zeta\Bigg].
\]
Applying the same test function:
\[
v_2(\zeta)=\eta_{z_{0}} \mathbf{1}_{R_{z_{0}}}(\zeta),
\]
this implies:
\[
\mathbb{E} \Bigg[\int_{R_{z_0}} \frac{\partial H_2}{\partial u_2}(\zeta)\eta_{z_{0}} d\zeta\Bigg]\leq 0.
\]
Since this holds for all \( z_0 \), we conclude:
\[
\frac{\partial H_2}{\partial u_2}(z_{0})=0.
\]
Thus, the first-order necessary conditions for the game problem are:
\[
\frac{\partial H_1}{\partial u_1} = 0, \quad \frac{\partial H_2}{\partial u_2} = 0.
\]
\fproof
\section{Applications: Pollution Control with Competing Regions}
Environmental regulation poses significant challenges when multiple regions share a common environment and must individually or collectively manage pollution levels. In many cases, regions act strategically, aiming to minimize their own costs, which include both the direct expenses of emission control and the environmental damage caused by pollution accumulation. The interaction between these regions creates a dynamic optimization problem that can be analyzed using stochastic control and game-theoretic approaches.

Pollution control often involves spatiotemporal dynamics, where emissions from one region affect the pollution levels in neighboring areas. This leads to complex models that incorporate both deterministic and stochastic effects, capturing environmental randomness such as fluctuating weather conditions, diffusion of pollutants, and uncertain effectiveness of control measures. In this section, we explore two different models that illustrate these challenges, starting from the following simple model.

\subsection{Example 1}

We consider a pollution control problem where two competing regions regulate emissions to minimize their respective costs. These costs include both the expense of controlling emissions and the environmental damage caused by accumulated pollution. Each region independently manages its emissions, leading to a strategic decision-making problem.

The evolution of the pollution level $Y(z)=Y(t,x)$ is modeled by the stochastic differential equation:
\begin{equation*}
    dY(z) = (u_1(z) + u_2(z)) dz + \sigma B(dz), \ \  t,x>0; \ \ Y(0,x)=Y(t,0)=y,
\end{equation*}
i.e.
\begin{align}
\frac{\partial Y(t,x)}{\partial t \partial x} =  u_1(t,x) + u_2(t,x) + \sigma \frac{\partial^2 B(t,x)}{\partial t \partial x},\ t,x>0; \  \  Y(0,x)=Y(t,0)=y,
\end{align}
where $Y(z)$ represents the pollution level at $z$, and $u_1(z), u_2(z)$ are the control efforts (e.g., emission reductions) by regions 1 and 2, respectively. The noise term represented by $\sigma B(dz)$ introduces stochastic fluctuations, accounting for unpredictable environmental factors such as weather conditions or industrial variability.

Each region seeks to minimize its expected cost, given by:
\begin{equation*}
    J_i(u_i) = \mathbb{E} \left[ \int_{R_Z} a_i u_i^2( \zeta) d\zeta+ c_i Y^2(Z) \right], \quad i=1,2.
\end{equation*}
The parameter $a_i > 0$ represents the cost of implementing pollution control measures, meaning higher values imply greater expense in reducing emissions. The term $c_i > 0$ reflects the sensitivity of region $i$ to pollution damage, where larger values indicate stronger concerns about environmental degradation.

To determine the optimal emission strategies, we introduce the Hamiltonian function:
\begin{equation*}
    H_i = a_i u_i^2 + (u_1 + u_2) p_i + q_i + (L_i * (u_1 + u_2))(z).
\end{equation*}
This formulation captures the trade-offs between emission control costs and pollution impact while incorporating the interactions between the two competing regions.

This problem follows a \textit{Linear-Quadratic} (LQ) framework. The pollution dynamics are linear in the control variables, while the cost function contains quadratic terms penalizing both emission reductions and pollution levels. The LQ structure ensures analytical tractability and provides insights into optimal control strategies for emission regulation.
Then, the first-order condition for optimal control $u_1$ is
\begin{equation*}
\frac{\partial H_1}{\partial u_1} = 0, \quad \textnormal{ which occurs when } \quad  2a_1 u_1 + p_1 = 0.
\end{equation*}
Solving for $u_1$, we obtain the optimal control
\begin{equation*}
\hat{u}_1 = - \frac{p_1}{2a_1}.
\end{equation*}
The adjoint variable $L_1$ is given by
\begin{equation*}
L_1 = -\frac{\partial H_1}{\partial y} = 0,
\end{equation*}
and the adjoint variable triple $(p_1,q_1,r_1)$ is the solution of the backward stochastic differential equation:
\begin{equation*}
\label{p1}
\begin{cases}
    p_1(dz)=- \frac{\partial H_1}{\partial y} dz + q_1(z) B(dz)
     + M_{1}^{r_1}(dz),\quad 0 \leq z=(t,x) \leq Z=(T,X), \\
    p_1(Z)=\frac{\partial g_1}{\partial y}(Y(Z)).
    \end{cases}
\end{equation*}
or
\begin{equation*}\label{p_1}
\begin{cases}
    p_1(dz)=q_1(z) B(dz)+ M_{1}^{r_1}(dz),\quad 0 \leq z=(t,x) \leq Z=(T,X), \\
    p_1(Z)=2c_1Y(Z).
\end{cases}
\end{equation*}
Similarly,  we get
\begin{equation*}
\frac{\partial H_2}{\partial u_2} = 0 \quad \textnormal{for } \quad  \hat{u}_2 =- \frac{p_2}{2a_2},\quad  L_2=0,
\end{equation*}
and
\begin{equation*}\label{p_2}
\begin{cases}
    p_2(dz)=q_2(z) B(dz)+ M_{2}^{r_2}(dz),\quad 0 \leq z=(t,x) \leq Z=(T,X), \\
    p_2(Z)=2c_2Y(Z).
\end{cases}
\end{equation*}
Comparing the equations for $u_1, L_1, p_1, q_1,$ and $r_1$, we see that
\begin{equation*}
p_2 = \frac{c_2}{c_1} p_1, \quad \textnormal{and} \quad\widehat{u}_2 =- \frac{p_2} {2a_2 } =- \frac{c_2 p_1}{c_1 2 a_2}= \frac{c_2 a_1}{c_1 a_2}\widehat{u}_1.
\end{equation*}
We conclude that a Nash equilibrium is obtained when the two regions coordinate their control efforts according to this proportion.

The original problem, formulated as a game between two regions controlling pollution emissions, reduces to a single-agent LQ control problem. Since the optimal control effort \( \widehat{u}_2 \) is expressed proportionally in terms of \( \widehat{u}_1 \), the system can be rewritten in terms of a single control variable, leading to a standard LQ optimization problem.

Substituting \(\widehat{u}_2\) into the pollution dynamics gives:
\begin{equation*}
    dY(z) = \left( u_1(z) + \frac{c_2 a_1}{c_1 a_2} u_1(z) \right) dz + \sigma B(dz),
\end{equation*}
which simplifies to:
\begin{equation*}
    dY(z) = \alpha u_1(z) dz + \sigma B(dz), \quad Y(0, x) = y,
\end{equation*}
where we define:
\begin{equation*}
    \alpha = 1 + \frac{c_2 a_1}{c_1 a_2}.
\end{equation*}
Summing the cost functions of both regions gives:
\begin{equation*}
    J(u_1) = \mathbb{E} \left[ \int_{R_Z} \left( a_1 u_1^2(\zeta) + a_2 \left(\frac{c_2 a_1}{c_1 a_2} u_1(\zeta) \right)^2 \right) d\zeta + (c_1 + c_2) Y^2(Z) \right].
\end{equation*}
Rewriting:
\begin{equation*}
    J(u_1) = \mathbb{E} \left[ \int_{R_Z} \left( a_1 u_1^2(\zeta) + \frac{c_2^2 a_1^2}{c_1^2 a_2} u_1^2(\zeta) \right) d\zeta + (c_1 + c_2) Y^2(Z) \right].
\end{equation*}
Defining:
\begin{equation*}
    \beta = a_1 + \frac{c_2^2 a_1^2}{c_1^2 a_2},
\end{equation*}
we obtain the simplified cost functional:
\begin{equation*}
    J(u_1) = \mathbb{E} \left[ \int_{R_Z} \beta u_1^2(\zeta) d\zeta + (c_1 + c_2) Y^2(Z) \right].
\end{equation*}
The Hamiltonian function is:
\begin{equation*}
    H = \beta u_1^2 + \alpha u_1 p + q.
\end{equation*}
The first-order condition for optimality is:
\begin{equation*}
    \frac{\partial H}{\partial u_1} = 2\beta u_1 + \alpha p = 0.
\end{equation*}
Solving for \({u}_1\):
\begin{equation}
{u}_1 = -\frac{\alpha p}{2\beta}.
\end{equation}
The adjoint equation follows from the costate dynamics:
\begin{equation*}
    dp(z) = q(z) B(dz) + M^{r}(dz), \quad p(Z) = 2(c_1 + c_2) Y(Z).
\end{equation*}
Thus, the optimal control policy is:
\begin{equation*}
    \widehat{u}_1 = -\frac{\alpha}{2\beta}  \widehat{p}.
\end{equation*}
Since \(\widehat{u}_2\) is expressed in terms of \(\widehat{u}_1\), this provides the full solution.
For further details, see Section 6.2 in Agram et al. \cite{AOPT1}. 
 We have proved:
\begin{theorem}
The optimal emission reductions $\widehat{u}$ 
\begin{equation*}
\widehat{u}_2 = \frac{c_2 a_1}{c_1 a_2}\widehat{u}_1, 
\end{equation*}
while \begin{equation}
    \widehat{u}_1 = -\frac{\alpha}{2\beta}   \widehat{p}.
\end{equation}
\end{theorem}

\subsection{Example 2}
We consider a two-region environmental control problem where each region aims to mitigate the effects of pollution over a shared geographical space \( \zeta \in R_Z \) and within a time horizon \( t \in [0, T] \). The pollution density \( Y(t,x) \) evolves according to a SPDE, which accounts for both deterministic dynamics (e.g., diffusion and control efforts) and random fluctuations (e.g., environmental noise). The problem is motivated by real-world scenarios in which multiple regions or entities must coordinate or compete in managing pollution levels, such as air quality management, environmental policies across borders, or shared water bodies affected by runoff.

Pollution is often a transboundary problem, where emissions from one region can affect the neighboring areas. Regions may implement various control strategies (e.g., emission reductions, investments in pollution abatement technologies, or policy measures), but such efforts are costly, and the effects are often uncertain due to environmental randomness (e.g., weather fluctuations, unforeseen industrial discharges). Hence, each region must strike a balance between reducing pollution and minimizing the cost of control efforts.

The evolution of the pollution density \( Y(t,x) \) at location \( x \) and time \( t \) is governed by the following SPDE:
\[
\frac{\partial Y(t,x)}{\partial t {\partial x}} =  - \alpha_1 u_1(t,x) - \alpha_2 u_2(t,x) + S(t,x) + \sigma(t,x) \frac{\partial^2 B(t,x)}{\partial t \partial x},
\]
where {$u_1(t,x)$ and $u_2(t,x)$ denote the control efforts by region 1 and 2, respectively, $a_1,a_2$ are positive constants, $(t,x)$ is the natural source of pollution, and $\sigma(t,x) \frac{\partial^2 B(t,x)}{\partial t \partial x}$ accounts for stochastic environmental noise.}

Each region \( i \) aims to minimize its individual cost function:
\[
J_i(u_1, u_2) = \mathbb{E} \left[ \int_{R_Z} \left( \frac{1}{2} Y^2(\zeta) + \frac{\beta_i}{2} u_i^2(\zeta) \right) d\zeta \right],
\]
where:
\begin{itemize}
    \item The term \( \frac{1}{2} Y^2(\zeta) \) penalizes {pollution levels}. The cost increases quadratically with the pollution density at location \( \zeta \), reflecting the negative impact of pollution on the environment, public health, and economic activity. For instance, high pollution may lead to increased healthcare costs, reduced quality of life, and environmental damage.
    
    \item The term \( \frac{\beta_i}{2} u_i^2(\zeta) \) penalizes the {cost of control efforts}. The parameter \( \beta_i \) determines how costly it is for region \( i \) to implement a particular level of control strategy \( u_i \). A higher \( \beta_i \) reflects more expensive control technologies or more stringent policies. The quadratic nature of this cost implies diminishing returns from control efforts, i.e., the more a region invests in control, the higher the associated costs.
\end{itemize}
This setup describes a two-region optimal control problem where each region \( i \) controls its own pollution level through the control \( u_i(t,x) \), but the two regions' actions affect the shared pollution dynamics. The challenge is to find the optimal control strategies \( u_1(t,x) \) and \( u_2(t,x) \) that minimize the expected cost \( J_1(u_1, u_2) \) and \( J_2(u_1, u_2) \), respectively.

The Hamiltonians for the two regions are given as:
\[
H_1 = \frac{1}{2} y^2 + \frac{\beta_1}{2} u_1^2 + p_1 \left(- \alpha_1 u_1 - \alpha_2 u_2 + s \right)+ \sigma q_1 +L_1 \star \left( - \alpha_1 u_1 - \alpha_2 u_2 + s \right),
\]
\[
H_2 = \frac{1}{2} y^2 + \frac{\beta_2}{2} u_2^2 + p_2 \left( - \alpha_1 u_1 - \alpha_2 u_2 + s \right)+ \sigma q_2 +L_2 \star \left( - \alpha_1 u_1 - \alpha_2 u_2 + s \right).
\]

The adjoint equations for \( p_1,q_1 \) and \( p_2,q_2 \) are derived from the state dynamics. For region 1, the adjoint equation is:
\begin{align}\label{p1a}
\begin{cases}
    p_1(dz)&=- Y(z)dz\\
    &+q_1(z) B(dz) 
     + M_{1}^{r_1}(dz),\quad 0 \leq z=(t,x) \leq Z=(T,X), \\
    p_1(Z)&=\frac{\partial g_1}{\partial y}(Y(Z)).
    \end{cases}
    \end{align}
and for region 2:
\begin{align}\label{p1a2}
\begin{cases}
    p_2(dz)&=- Y(z)dz\\
    &+q_2(z) B(dz) 
     + M_{2}^{r_2}(dz),\quad 0 \leq z=(t,x) \leq Z=(T,X), \\
    p_2(Z)&=\frac{\partial g_2}{\partial y}(Y(Z)).
    \end{cases}
    \end{align}
The optimality condition for \( u_1 \) is:
$$\frac{\partial H_1}{\partial u} = \beta_1 u_1- \alpha_1 p_1-(L_1 \star \alpha_1)=0.$$
Hence, we must have
$$
u_1(z)=\frac{\alpha_1 p_1-(L_1 \star \alpha_1)}{\beta_1}.
$$
Since
\begin{align*}
    \frac{\partial H_1}{\partial y}= \frac{\partial H_2}{\partial y}=y.
 \end{align*} 
The adjoint variables \( L_1=L_2 \) are explicitly given by:
\begin{align}\label{leq}
L_1(z) =L_2 (Z)= -Y(z).
\end{align}
Using this, \( (L_1 \star \alpha_1)(z) \) simplifies to:
\[
(L_1 \star \alpha_1)(z) = \alpha_1 x(T-t) L_1(z) + \alpha_1 \int_0^x \int_t^T L_1(\zeta') d\zeta'.
\]
Substituting \( L_1 \) into the expression for \( u_1 \):
\[
u_1 = \frac{\alpha_1 p_1 - \alpha_1 x(T-t)L_1(z) - \alpha_1 \int_0^x \int_t^T L_1(\zeta') d\zeta'}{\beta_1}.
\]
Similarly, for \( u_2 \):
\[
\frac{\partial H_2}{\partial u_2} = \beta_2 u_2 - \alpha_2 p_2-(L_2 \star \alpha_2) = 0 \quad \Rightarrow \quad u_2 = \frac{\alpha_1 p_1-(L_2 \star \alpha_2)}{\beta_2},
\]
where \[
(L_2 \star \alpha_2)(z) = \alpha_1 x(T-t) L_1(z) + \alpha_1 \int_0^x \int_t^T L_1(\zeta') d\zeta'.
\]
Thus
\[
u_2(z) = \frac{\alpha_2 p_2 - \alpha_2 x(T-t)L_2(z) - \alpha_2 \int_0^x \int_t^T L_2(\zeta') d\zeta'}{\beta_2}.
\]
Both controls \( u_1 \) and \( u_2 \) reflect the interplay between their respective adjoint variables, \( p_1 \) and \( p_2 \), and the pollution dynamics, including diffusion and stochastic effects. The optimal strategies balance the cost of control efforts with the effectiveness of reducing pollution in the shared domain.

We summarise this as follows:
\begin{theorem}
Assuming that $ \beta_1\neq 0$ and $ \beta_2\neq 0$. Then the optimal control values for each region $(\widehat{u}_1,\widehat{u}_2)$ are given by the following expressions, respectively:
\[
\widehat{u}_1(z) = \frac{\alpha_1 p_1 - \alpha_1 x(T-t)L_1(z) - \alpha_1 \int_0^x \int_t^T L_1(\zeta') d\zeta'}{\beta_1}.
\]
and 
\[
\widehat{u}_2 (z)= \frac{\alpha_2 p_2 - \alpha_2 x(T-t)L_2(z) - \alpha_2 \int_0^x \int_t^T L_2(\zeta') d\zeta'}{\beta_2}.
\]
where the adjoint variables $p_i$ and $L_i$, for $i=1,2$, satisfy the equations \eqref{p1a}-\eqref{p1a2} and \eqref{leq}, respectively.
\end{theorem}

\subsubsection{Symmetry in Parameters}

In the special case where \(\alpha_1 = \alpha_2 = \beta_1 = \beta_2 = 1\), the problem exhibits perfect symmetry between the two players. The optimal controls simplify to:
\[
u_1 (z)= p_1 - x(T-t) L_1(z) - \int_0^x \int_t^T L_1(\zeta') d\zeta',
\]
\[
u_2 (z)= p_2 - x(T-t) L_2(z) - \int_0^x \int_t^T L_2(\zeta') d\zeta',
\]
where \(L_1(z) = L_2(z) =-Y(z)\).

This symmetry indicates that the two regions exert equal influence on pollution dynamics. The controls reflect a trade-off between immediate effort and the anticipated benefit of reducing pollution. Given identical parameters, the behavior of the two regions differs only due to their adjoint variables \(p_1\) and \(p_2\), which are shaped by initial conditions, external pollution sources, and stochastic fluctuations.

Implications for players:
\begin{itemize}
\item Equal roles: Both regions have identical strategies and impacts on pollution mitigation.
\item Balanced efforts: The control efforts balance cost and pollution reduction uniformly across regions.
\item Outcomes depend on external Factors: Differences in outcomes arise only from initial pollution levels, external sources, or environmental noise, as the inherent dynamics are otherwise identical.
\end{itemize}
This case highlights the interplay of stochastic dynamics and control in scenarios where multiple agents with symmetric roles aim to mitigate shared environmental challenges.

\[\]
\textbf{Acknowledgments.}
Nacira Agram  gratefully acknowledge the financial support provided by the Swedish Research Council grants (2020-04697), the Slovenian Research and
Innovation Agency, research core funding No.P1-0448. Olena Tymoshenko acknowledges support from the MSCA4Ukraine project (AvH ID:1233636), which is funded by the European Union. Views and opinions expressed are however those of the author only and do not necessarily reflect those of the European Union. Neither the European Union nor the MSCA4Ukraine Consortium as a whole nor any individual member institutions of the MSCA4Ukraine Consortium can be held responsible for them.

\end{document}